# Increasing student engagement in math placement and preparation

Debra Lewis

Mathematics Department, University of California Santa Cruz, Santa Cruz, CA 95064, USA
Email: lewis@ucsc.edu

Abstract

Math placement is a crucial step between admission and full engagement with the university for students in STEM majors. Students' placement experiences influence not only their mathematical entry point on their degree pathway, but their perceptions of their relationship with the university. Application of behavioral economics theory, particularly choice architecture, to placement communication can increase student engagement with the placement process, reducing enrollment in preparatory courses while maintaining positive outcomes in more advanced courses. Targeted communication with students facilitates successful placement and enhances self-efficacy and sense of belonging. Math placement is a pigeonholing process and inherently involves an array of sub-instructions and information relevant only to subgroups; optimizing content flow improves retention of key information and minimizes exposure to intimidating or alienating information.

Keywords: Math placement, choice architecture, behavioral economics, targeted communication, reassessment.

## 1. Introduction

For students intending to major in STEM disciplines, math placement is a critical step between admission and full engagement with the university. Accurate, timely placement and enrollment provides students with an entry point to their STEM development that facilitates rapid, productive progress towards degree. Clear, supportive communication with students not only facilitates successful placement, but provides substantive evidence of the university's commitment to students; for better or worse, the placement experience influences students' academic expectations and perceptions of college.

The predictive validity of placement assessments is significantly reduced if students have insufficient knowledge about the placement process (Hodara, Jaggars & Karp, 2012; Briggs, 2009). Misperceptions about the purpose of the placement process, inadequate knowledge about preparation for placement, and low math self-efficacy deter many students from preparing for their placement assessment (Fay, Bickerstaff & Hodara, 2013; Hoffman, 2015; Venezia, Bracco & Nodine, 2010). Many students fear that reviewing prior to placement assessment could result in placement in courses beyond their academic abilities, and hence avoid even minimal review of previously mastered material (Jaggars & Hodara, 2011).

The first step in increasing the accuracy and equity of the placement process is to acknowledge that any such process is imperfect, and that performance in math assessment is strongly influenced by non-cognitive factors, some of which can be can be exploited, others that can be ameliorated (e.g. triggers for stereotype vigilance and math anxiety (Spencer, Steele & Quinn, 1999; Foley et al., 2017)), and others over which the university has no control (e.g. an all-around lousy day). It is then possible to evaluate and design assessment tools, protocols, and communication with an eye to maximizing the likelihood of accurate placement and minimizing the consequences of a suboptimal placement experience.

Communication with incoming students in the context of math placement offers an underutilized opportunity to positively influence student expectations and perceptions of their college experience. Placement support should not only transmit essential information about the process, but should stimulate and productively channel incoming students' enthusiasm about their college education, communicate the university's appreciation of their academic potential, boost their math self-efficacy, and increase their willingness to engage with unfamiliar pedagogic modalities. The instructional structures central to most "first year" university math courses—very large lectures, discussion sections, and office hours—are new and potentially alienating to most incoming students; a positive experience with the math placement process can build student trust and sense of belonging, which in turn benefits academic performance and persistence (Kuh et al., 2008; Lent, Brown & Larkin, 1986; Multon, Brown & Lent, 1991).

*1.1 The costs of placement or enrollment below true skill and knowledge level*

While the risks of placing students into too advanced a course are widely recognized, placing students below their true skill and knowledge level is also dangerous. Overestimating a student's current knowledge and skills can set them up for failure in their first math course; underestimating a student's current abilities can retard academic progress, undermine self-efficacy, and lead to boredom and academic disengagement. The cumulative cost to students and institutions of placement or enrollment below true capabilities is enormous; in one urban community college system, it was estimated that one quarter of students who were placed into developmental math could have succeeded in college-level math (Scott-Clayton, 2012).

Enrollment below true experience and ability level disadvantages students in several ways:
- Students who are accurately placed into preparatory math instruction are vulnerable to being perceived as—or perceiving themselves as—"bad at math" because they need to work to learn material that many of their classmates already know.
- Finite resources are spread across larger groups of students, so students who could best benefit from additional support don't get as much as they could if preparatory courses had fewer advanced students.
- Students who deliberately spend a full academic term reviewing material they already know can be overwhelmed by their first STEM course (if any) in which a substantial percentage of the course content is unfamiliar.
- Students who plan to major in a STEM field risk an increased time to degree if they must complete preparatory coursework before beginning a calculus sequence.
- Students with low math self-efficacy who don't trust—or are told that they shouldn't trust—their placement or progress because of external factors (e.g. stereotypes about ethnicity and/or gender and mathematics) may be better served by instruction directly addressing those factors than by 'cautious' enrollment recommendations that validate their self-doubt.

Accurate placement, effective preparation and instruction, and supportive communication are crucial for educationally or economically disadvantaged students. Many students from underserved high schools have gaps in their basic mathematical preparation, but their mathematical self-efficacy and motivation can be further eroded if they are advised to repeat a course taken years earlier, particularly if that advice is not aligned with their placement, or the recommended coursework fails to address the pertinent deficiencies. Many students, particularly those from underrepresented or underserved groups, are advised to "play it safe" and enroll in preparatory courses, even though they've successfully completed the analogous high school course in high school and demonstrated mastery of the course content in placement; this supposedly cautious strategy can result in increased pressure to succeed in all subsequent math courses to avoid delayed major declaration or graduation, and can increase STEM attrition.

*1.2 Maximizing the benefits of innovation in placement*

Adaptive learning and assessment systems, including ALEKS PPL and Smarter Balanced, provide detailed estimates of student mastery of a broad range of concepts and skills (Falmagne & Doignon, 2011; McGraw Hill, 2019; Smarter Balanced Assessment Consortium, 2018). ALEKS PPL offers a blend of formative and cumulative assessment: the results of cumulative assessments are used both for placement and for construction of a personalized learning path for guided review, if needed; students can reassess (possibly multiple times) after a mandatory minimum amount of review. Students receive detailed feedback after each practice problem; retention of (re)learned content is monitored by means of frequent, brief assessments; once the mandatory review period has been completed, the student can initiate another placement assessment.

Review and reassessment are particularly valuable for students with previous math instruction of inconsistent rigor and efficacy (e.g., students from underserved schools). By first identifying gaps and weaknesses, and then providing students with the opportunity to repair those deficiencies pre-arrival, reassessment cycles enable many students who experienced windows of marginal K-12 math instruction to begin major-required mathematical coursework in their first term. Targeted self-study offers a feasible approach to (re)acquisition of elementary skills not covered in university-level courses. For example, directing a student who assesses as calculus-eligible, but stumbles when adding fractions, to a precalculus course is unlikely to be beneficial, given that addition of fractions isn't covered in precalculus courses.

To fully benefit from a placement process involving both assessment and guided improvement, students must understand, trust in, and exploit the differences between such a process and conventional standardized tests (Venezia, Bracco & Nodine, 2010). Adaptive placement systems, particularly those supporting review and reassessment, can be powerful tools for the development and reinforcement of mathematical self-efficacy and self-regulated learning (Hattie & Timperly, 2007; Nicol & Macfarlane-Dick, 2006; Sadler, 1998). However, their influence can be hobbled if external communication about the placement process undercuts and casts doubt on the validity of gains made via self-study and reassessment.

Assertions that a student's initial placement assessment most accurately reflects their 'true' abilities and knowledge, and that improvements resulting from review and reassessment are deceptive, may be manifestations of latent fixed mindset beliefs (Dweck & Leggett, 1988; Dweck, 2015). Characterization of enrollment according to improved placement as 'skipping' or 'jumping over' the lower level courses challenge the legitimacy of the placement process; uncertainty about the accuracy or legitimacy of assessment can reduce the benefits of performance-approach goals (Darnon et al., 2007). Skepticism regarding review and reassessment opportunities may reflect a fundamental confusion between or conflation of knowledge/skills assessment and 'aptitude' tests; few instructors would reject the results of a final exam as a measure of student readiness for the next course on the grounds that students had studied for the final.

*1.3 Non-cognitive challenges and benefits in math placement*

Low math self-efficacy and negative expectations regarding assessment can be strong deterrents to full engagement with the placement process. Accurate matching of students with the information pertinent to their situation enhances a sense of belonging and avoids the detrimental effects of generic statements on academic performance (Park et al., 2017); delivery of content targeting specific subgroups reduces the risks that essential information will be overlooked, and communicates the institutions commitment to that subgroup and awareness of its particular needs (Hodara, Jaggars & Karp, 2012).

Demographic factors and personal experiences play significant roles in math placement and subsequent course success. Stereotype threat in high stakes assessments is a serious burden for many under-represented minority and First Generation students (Maloney, Schaeffer, & Beilock, 2013; Steele & Aronson, 1995; Spencer, Steele & Quinn, 1999). Performance under pressure is, of course, a valuable skill in many contexts, but the pressure is far from equal for all students—societal expectations and personal/financial consequences vary dramatically. Math anxiety and stereotype threat have a significant negative effect on many students' performance in high stakes assessments; students from underrepresented minorities and underserved schools are particularly vulnerable to stereotype threat and 'inherited' math anxiety (Elliot & Thrash, 2004; Lea, 2010; Lyons & Beilock, 2011; Maloney, Schaeffer, & Beilock, 2013; Maloney et al., 2015; Sprute & Beilock, 2016).

Strong mathematical self-efficacy is crucial for academic success in STEM (Pajares, 1996, 2001; Shunk & Pajares, 2002; Jackson, 2002). The math placement process should enhance students' mathematical self-efficacy; communication about that process should alleviate uncertainty and, to the extent possible, minimize the stress of assessment and enrollment (Putwain, Woods & Symes, 2010). Stress impairs working memory and cognitive flexibility, and generally impairs long term memory retrieval (Shields, Sazma & Yonelinas, 2016); by increasing student stress levels, intimidating communication about math placement can degrade cognitive function and self-efficacy in vulnerable students (Maloney, Schaeffer & Beilock, 2013; Maloney, Sattizahn & Beilock, 2014; Herts & Beilock, 2017; Zajacova, Lynch & Espenshade, 2005).

The cognitive and physiological effects of math anxiety may be strongest when anticipating a high stakes math task, not when executing one (Lyons & Beilock 2012); thus the challenge is to successfully transmit key information when anxious students' capacity to process and act on that information is impaired. Students who feel frustrated, neglected, out of place, or incompetent because they can't find the answers to their questions, need to wait several days for assistance, or are reproached for asking questions face an increased risk of suboptimal performance due to stress and reduced sense of belonging. Stress degrades working and long-term memory, impairing performance in assessment (Ashcraft & Ridley, 2005; Ashcraft & Krause, 2007; Shields, Sazma & Yonelinas, 2016; Chang & Beilock, 2016; Herts & Beilock, 2017). Insight into the impact of avoidance, persistence, and informed commitment to academic goals, e.g., reluctance to reassess or enroll in a math course stemming from uncertainty about the placement process or courses, can guide improvements communication about requirements and expectations.

## 2. Material and methods

The University of California at Santa Cruz (UCSC) adopted the math placement and preparation system ALEKS PPL in 2015 on the strength of the accuracy of placement over a wide range of skill levels and the integrated self-study review and learning support. The institutional practices and communication tools for math placement used at UCSC prior to the adoption of ALEKS PPL were developed to support a single placement assessment during a narrow window, and were poorly matched to the new placement process. The development and implementation of a portfolio of placement communication protocols and tools encouraging full student engagement with ALEKS provided an opportunity to study the influence of communication modalities and choice architecture on placement activity and outcomes.

We first provide a brief overview of the academic context and math placement process at UC Santa Cruz, We then describe the communication strategies that were implemented with the goal of increasing productive engagement with the process and shifting enrollment from preparatory courses (*College Algebra* and *Precalculus*) to calculus courses, and the methods used to assess the efficacy of the new communication resources and protocols.

### 2.1 Math placement at UC Santa Cruz

UCSC is one of ten campuses in the University of California system. UCSC is a Hispanic Serving Institution with a diverse student population. 37% of UCSC undergraduates are first generation college students, 37% are Pell grant-eligible, and 32% are underrepresented minorities (University of California, 2019). Many domestic students come from affluent Bay Area communities, while others come from rural agricultural communities in the Central Coast region. Roughly half of UCSC's undergraduates major in STEM disciplines.

UCSC's math enrollment policies implement a unidirectional guided self-placement process: students have *downward* freedom of choice, but no upward flexibility. Students can enroll below their placement—but not above it—and students who initially place into a developmental course are not required to review and reassess prior to enrollment. Student surveys, with self-reported high school coursework, indicate that a majority of incoming UCSC students who take a mathematics course in their first term repeat coursework already completed in high school. (See Figure A.1, Appendix A.)

Differences in previous educational opportunities is a crucial factor in the disparities in placement across ethnic and socio-economic groups. The broad range of socio-economic backgrounds for UCSC students result in widely varying K-12 experiences, with white and Asian students being less likely than Hispanic students to come from underserved school districts. Students from rural agricultural areas frequently have fewer options in their high school math courses than their peers from wealthier districts in tech-centric areas. Hence the efficacy of a review-and-reassessment placement process in narrowing equity gaps resulting from differences in prior preparation is of particular relevance at UCSC.

An overwhelming majority of UCSC students who review and reassess in ALEKS PPL improve their enrollment opportunities. Six out of seven UCSC students who initially place into *College Algebra* improve their enrollment eligibility if they review and reassess, despite the fact that the lowest placement tier spans three fifths of the score range. Nearly two thirds of the reassessing students from this tier become calculus-eligible. Twenty out of twentyone reassessing students who initially place into *Precalculus* become calculus-eligible.

Once demographic characteristics–race/ethnicity, gender, first generation status and Pell grant recipient status–and academic preparation (as measured by SAT Math scores and high school GPA) are taken into account, students who become eligible to enroll in a more advanced course after reassessing in ALEKS PPL are as likely to pass, and as likely to earn a B or better, as students who placed into that course with their initial placement assessment (Fernald, Lewis & Padgett, 2016). No differences in grades in the subsequent math class between students who reassessed and those placed directly into their first course have been found.

## 2.2 Engagement enhancement strategies—marketing placement investment

UCSC's communication with students regarding math placement traditionally had an authoritarian, at times coercive, tone. Imperative formulations and emphasis on negative repercussions of non-compliance and poor performance were pervasive. Information was organized from an administrative perspective, which often failed to align with student needs/expectations, making it difficult for students to find, filter, prioritize, and retain information. Use of specialized terminology in student-facing webspace and documents resulted in misunderstandings that reduced student confidence in UCSC or themselves.

Two years of UCSC placement activity and outcomes data, supplemented by student surveys and interviews, guided the development in 2017 of an online placement support resource, the UCSC Math Coach. Prioritization and organization of placement and enrollment information and advice was based on student feedback and queries, yielding as-needed guidance and answers to common questions. Student interviews provided essential guidance for the presentation of essential reminders and information about reassessment outcomes.

### 2.2.1 Goal-orientation and choice architecture

Nudge theory, or libertarian paternalism, focuses on non-coercive influences on human behavior (Sunstein, 2014; Ebert & Freibichler, 2017). "A nudge, as we will use the term, is any aspect of the choice architecture that alters people's behavior in a predictable way without forbidding any options or significantly changing their economic incentives... Putting fruit at eye level counts as a nudge. Banning junk food does not." (Thaler & Sunstein, 2008). Three key nudge techniques that increase the likelihood that individuals will choose the desired options are defaults, social proof heuristics, and increased salience. These techniques can be particularly powerful in contexts in which individuals have little knowledge or understanding of the possible relative costs and benefits of the available options, and weak pre-existing biases regarding the alternatives.

Choice architecture can be exploited in placement communication to increase student engagement and utilization of available resources by increasing the salience of the most promising options and exploiting the instinctive belief in the safety of following the herd. For students with low self-efficacy or stereotype vulnerability, social proof heuristics may be more protective than persuasive; a decision consistent with an assertion of mathematical competence can be justified as "the usual thing" even if self-efficacy is relatively low or the decision might otherwise be interpreted as defying prevailing stereotypes (with the associated heightened penalty for suboptimal performance). Opt in/opt out nudges can be highly effective when the relevant decision has negative associations for many members of the target population; e.g., considering registering as an organ donor invokes the possibility of an early, abrupt death, so changing from opt in to opt out enrollment in a donor program can dramatically increase the number of potential donors (Thaler and Sunstein, 2008).

Many colleges and universities show students a math placement conversion, or look up, table, with score thresholds were listed in the leftmost column and course eligibility listed to the right. (See Figure B.1, Appendix B.) Such a table is naturally interpreted as saying "If your score is in this range, then you can enroll in this/these courses." A starting-point oriented placement score conversion table embodies a choice architecture that *discourages* reassessment by subordinating students' academic goals to their immediate circumstances. Reversing the order of presentation, so as to say "If this is the first math course required for your major, then do your best to earn a score in this range" focuses students' attention on their long term goal—completion of a degree in their chosen field—and presents placement preparation and assessment as a means to that end. (See Appendix B, Figure B.2.)

Goal-oriented placement communication supports a growth mindset and student engagement (Horstmanshof & Zimitat, 2007). In the context of a review-and-reassess placement process, goal-oriented nudges convey institutional confidence in students' capacity to regain and reinforce basic skills and knowledge outside of the conventional classroom setting.

*2.2.2 Students as agents of the university—establishing the contract*

Economic agency theory addresses problems arising when an agent acts for a principal: the preferences or goals of the principal and agent may differ or be in direct conflict; the principal and agent may have different perceptions of the risks involved; if may be difficult or expensive for the principal to gain information about the agent's actions (Jensen & Meckling, 1976; Eisenhardt, 1989). If we regard the student as the university's agent in determining the appropriate first math course for the student, these problems are clearly inherent to the placement process. Economic agency theory offers valuable insights into the administrative and psychological challenges of implementation of an efficient, effective placement process (but should not be conflated with theories of student agency).

While the student and the institution hopefully share long range, high level goals—academic success leading to personal growth and a fulfilling career—short term goals and prioritization of risks can differ dramatically. Both individuals and institutions demonstrate delay discounting, but the time scales involved may be dramatically different; among students, socioeconomic factors and career goals strongly influence perceptions of risk and opportunity. Interactive learning systems such as ALEKS PPL provide detailed information about student activity in the system, but it can be challenging to extract actionable insights from that plethora of data, and little or no information is available about students' preparatory activities outside of institution-linked online resources.

Presenting the placement process as a service or benefit that the university offers to the student (accurate placement increases the likelihood that the student will make substantials gains in skill and knowledge , and will progress towards their academic goals) is unquestionably preferable to presenting placement as a constraint unilaterally imposed on the student, but it is advantageous to bear in mind that the benefit comes at a cost to the student—the student earns valuable information about alignment between their skills and the university's courses by undertaking a stressful and intellectually demanding activity. Asking students to review and reassess if their initial assessment score is low significantly increases the immediate workload; hence the need to clearly present the likely compensation for this effort. Incentives can increase agency and reduce irrational behavior driven by avoidance of psychological pain (Jensen, 1994). Thus compelling incentives can potentially reduce the negative consequences of the pain experienced by highly math anxious students when anticipating a math assessment (Lyons & Beilock, 2013).

An example of proscriptive phrasing in the placement contract involves the use of 'only if' statements rather than 'if' statements. A commitment to ensure a desired outcome *if* specified conditions are met is a different factual statement, and presents a different balance of social power between the student and the university, from an admonishment that the outcome can be obtained *only if* those conditions are met.

As an illustrative example, we consider an example from advising email communication with incoming students about math placement. Incoming students have an enrollment window during their orientation visit (students are told this in other email); after that window closes, they are unable to enroll in courses for approximately two weeks, but enrollment opens again approximately six weeks before the beginning of the term and remains open through the first few weeks of the fall term (this is not stated in advising email); students who improve their placement after their orientation session can change their math enrollment. Compare the 'if' statement

> *If you complete the placement process by [date], you will be able to enroll in your math course during Summer Orientation.*

to the 'only if' statement

*You must complete the placement process by [date] to enroll in a math course.*

The 'if' version offers an incentive—by taking timely action, the student can complete an important task and avoid uncertainty about course availability—and tacitly conveys the important information that it is possible to continue the placement process beyond that target date. Explicitly linking the soft target date for completion of math placement to Summer Orientation increases transparency of the process. The 'only if' version introduces a spurious hard deadline and penalty (inability to enroll in a math course for the fall term) for failure to meet that deadline, and offers no context or justification for that policy.

Adoption of an online assessment and study system supports broad flexibility in timing of assessment; this flexibility introduces new opportunities and challenges. Last minute assessment is rampant, but not irrational, and should not be dismissed as procrastination or poor time management. In economics, time preference and delay discounting model the decrease in perceived value of goods or benefits valuation received at a later date (Doyle, 2013). If you would rather be given $5 dollars now than $10 dollars a month from now, you've shown a strong preference for an immediate payoff—or are coping with extreme inflation. While the strength of delay discounting is, to some extent, a personal characteristic, external factors influence the degree of discounting; reduced executive abilities limit the individual's ability to resist the tendency to favor the present (Bayer & Osher 2018).

Delay discounting may play a crucial role in many students' limited engagement with the placement process. Pragmatic, purely rational optimization suggests that a math-anxious student who knows that the odds are very high that they can satisfy a math requirement by completing an assessment sufficiently far in advance of a deadline that they can review for several hours and then reassess would do that, rather than spending an entire academic term in preparatory course; students with little confidence in their mathematical abilities should, on purely logical grounds, prefer a non-credit bearing preparatory program to a several-credit course that will be part of their university transcript. However, the last few days before assessments must be completed for use in priority enrollment during Summer Orientation are by far the busiest of the year. The central temporal consideration in math placement for incoming frosh is the constraint that students must satisfy course prerequisites prior to enrollment.

Intimidating communication regarding placement timelines can backfire, *increasing* assessment avoidance among math-anxious students by increasing their stress levels and focusing their attention on the negative consequences of possible poor performance. In individuals who identify as math-anxious, anticipation of a math question has been shown to trigger activity in the part of the brain that processes physical pain (Lyons & Beilock, 2012). This suggests that delayed assessment may not be due to poor time management, insufficient motivation, or failure to understand instructions, but to a combination of an entirely rational desire to avoid pain and anticipation of a poor outcome.

### 2.2.2 Targeted information flow

Math placement is a pigeonholing process, and inherently involves an array of sub-instructions and information relevant only to some students. Exposure to content relevant only other student groups can dilute the impact of pertinent information; apparent prioritization of 'others' can unintentionally convey that the reader is out of place or undervalued. For example, while students who have already satisfied placement requirements via AP or IB test credit do not need to know about alternative means for completion of placement (and—if there is a cost per student for that alternative— the university may wish to deter superfluous assessment), detailed presentation early in a one-size-fits-all placement communication of information relevant only to AP/IB test takers may signal to other students that their mathematical preparation is deficient, leading to increased anxiety about comparison to the 'others' in courses, and subsequent avoidance of those courses. Generic statements about coarse student groupings can impair academic performance (Park et al., 2017).

Effective use of navigation within websites and adaptive email can limit exposure to irrelevant content and provide students with an intuitive progression through relevant content. For example, if a review and reassess placement process is used, it can be advantageous to bifurcate post-assessment communication into not-yet-calculus-eligible and calculus-eligible information streams, presenting information about self-remediation and likelihood of improving enrollment eligibility to students placing into preparatory course-only eligibility, and presenting information about entirely optional, but advantageous, preparation for peak performance to calculus-eligible students.

Respectful, informative responses to student requests for clarification or assistance are essential. If students feel neglected, or reprimanded for failure to understand or correctly follow the placement process, they may become angry and resentful, or may decide that they must not belong because they didn't see what they believe everyone else understands. Efficient placement problem resolution requires in-depth knowledge of both the way things actually work and the ways many students expect things to work. Avoiding reproachful and patronizing language, apologizing for confusing or incomplete information, and explaining the reasons for any counter-intuitive policies and inefficient processes, can make problem resolution a positive experience in which the student learns that they can overcome setbacks and successfully navigate the complexities of university bureaucracy.

## 2.3 Data collection and analysis

Student activity in ALEKS PPL was tracked using ALEKS comprehensive cohort reports, which include date and duration of placement assessments, placement scores, eleven subtopic scores, and—for students who review in their ALEKS PPL Learning Module and reassess, data tracking Learning Module activity. Course outcomes and ethnicity data were extracted from the UCSC Academic Information System.

An in-depth analysis of early outcomes for UCSC students who completed math placement using ALEKS PPL was performed by UCSC Institutional Research and Policy Studies (Fernald, Lewis & Padgett, 2016) as part of the University of California BFI Adaptive Learning Technology Pilot (Office of the President of the University of California, Institutional Research and Academic Planning, 2016).

Self-reported final high school math course data was obtained from an online survey of UCSC students who had taken a first year math course (*College Algebra, Precalculus, Calculus for the Social Sciences, Calculus with Applications, Calculus for Mathematics, Physics & Engineering,* or *Honors Calculus*) in Fall 2016. The same survey gathered information about student perceptions of the accuracy of their placement; preparation for placement (time spent, resources used); value, completeness, and accuracy of placement communication; and utilization and perceived value of instructional resources in their Fall math course (e.g., lectures, discussion sections, office hours, LMS). The survey was administered at the end of the Fall term, but prior to final exams, to capture students' self-estimates of their mastery of course content. 997 students completed the survey, slightly more than 40% of the students contacted.

In Winter 2017, an in-depth survey of one dozen first year students enrolled in an experimental adaptive offering of *College Algebra* solicited information about perceptions regarding math placement, particularly placement communication. These students were selected because most had failed *College Algebra* the previous term, and none had completed more than the mandatory initial placement assessment the previous summer; all made good-to-excellent progress in ALEKS in the context of the winter *College Algebra* course. Thus, they were considered to be students who could have benefited significantly from summer review and reassessment in ALEKS, but had been unaware of the opportunity (did not read or remember reading multiple email about reassessment) or skeptical that they would benefit from it. The survey was developed by the instructor of the adaptive *College Algebra* course and a UCSC Computer Science PhD candidate in the Assistive Sociotechnical Solutions for Individuals with Special needs using Technology (ASSIST) Lab.

## 3. Results and discussion

### *3.1 Shifts in placement and enrollment*

The number of students initially placing into the preparatory tiers decreased after the new placement resource, the UCSC Math Coach, was introduced, contributing to overall improvements in calculus-readiness. In the first two years after UCSC's adoption of ALEKS PPL, 65% of the students who assessed in ALEKS PPL were calculus-eligible upon completion of the placement process; in the two years since the implementation of the Math Coach, 87% were calculus-eligible upon completion of the placement process. Eligibility for the calculus sequence for engineering, physics, and mathematics increased from 55% to 69%.

Enrollment in the preparatory courses *College Algebra* and Precalculus has monotonically declined over the four years that UCSC has used ALEKS PPL for math placement, while fall term calculus enrollment increased dramatically with the adoption of ALEKS PPL and has remained approximately constant, as a percentage of the incoming class, since then. *College Algebra* enrollments in Fall 2018 were less than a third of those in Fall 2014, while Fall 2018 *Precalculus* enrollments were less than one third of those in Fall 2014; reductions were not imposed using enrollment caps—preparatory courses have not exceeded capacity since the adoption of ALEKS PPL.

|  | Before Coach: 2015 and 2016 | | | With Coach: 2017 and 2018 | | |
|---|---|---|---|---|---|---|
| Placement tier | first | best | net shift | first | best | net shift |
| 100 | 33% | 21% | -36% | 28% | 14% | -50% |
| 200 | 18% | 14% | -22% | 14% | 9% | -36% |
| 300 | 9% | 10% | 11% | 8% | 8% | 0% |
| 400 | 26% | 36% | 38% | 28% | 39% | 39% |
| 500 | 15% | 19% | 27% | 22% | 30% | 36% |

*Table 1. Placement tier results: percentage of students in each tier after initial assessment, after completion of placement (including students who assessed only once), and net shift from initial to best tiers. Enrollment eligibility caps: 100: College Algebra, 200: Precalculus, 300: Calculus with Applications (life sciences), 400: Calculus for Math, Physics, and Engineering, 500: Honors Calculus.*

Disaggregation by self-reported last high school math course suggests that reassessment captures typical retained/regained knowledge and skills from high school coursework. Assessment scores aggregated by final high school math course showed only slight differences between mean scores for students who completed a single placement assessment and the best assessment scores of students who completed multiple assessments (Appendix A). The mean differences between initial and best scores for students completing multiple placement assessments were far larger, ranging from an average improvement of 25 points (out of 100) for students whose reported last high school math course was Advanced Algebra, to 17 points for students whose last math course was (non-AP) calculus. This suggests that some students who reassess may have estimated what they "should" be able to do, based on their prior coursework, and reassessed because their initial score fell below that threshold.

| Academic year | Assessed using ALEKS PPL | Reviewed & reassessed | College Algebra | Fall Precalculus | Fall Calculus* |
|---|---|---|---|---|---|
| 2014-15 | pre-ALEKS | — | 13% | 26% | 21% |
| 2015-16 | 86% | 28% | 7% | 19% | 34% |
| 2016-17 | 74% | 24% | 6% | 11% | 33% |
| 2017-18 | 72% | 25% | 4% | 10% | 35% |
| 2018-19 | 67% | 20% | 3% | 8% | 34% |

*Table 2. Placement activity and enrollment for first year students, as percentages of the total number of these students. College Algebra is offered only in fall; Precalculus and Calculus are offered each quarter. *Calculus enrollment shown is only for the first courses in each of the two main calculus sequences: Calculus with Applications (life sciences), and Calculus for Mathematics, Physics, and Engineering.*

## 3.2 Enrollment below placement

While enrollments in preparatory courses have decreased significantly, placement distributions *within* those courses remain largely unchanged. Less than a third of the UCSC students who took *College Algebra* in 2016-18 reassessed before taking the course, and some students who did reassess into a higher placement tier remained enrolled in *College Algebra* despite their eligibility to enroll in a more advanced course.

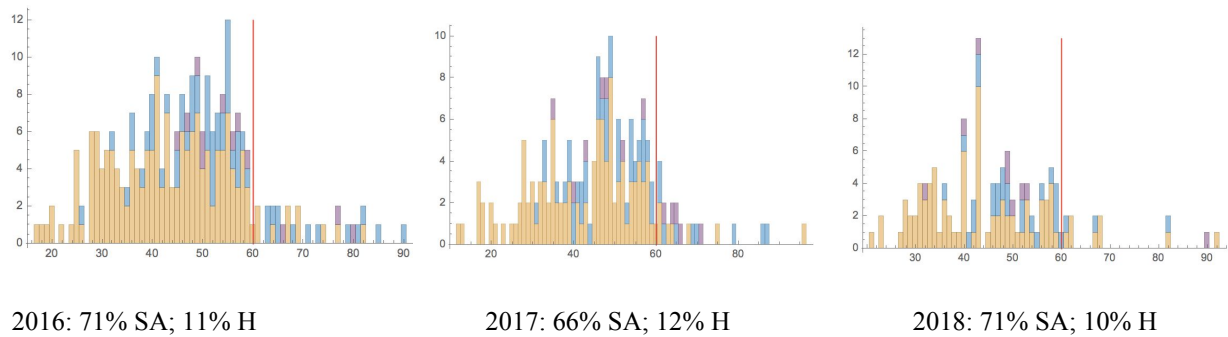

| 2016: 71% SA; 11% H | 2017: 66% SA; 12% H | 2018: 71% SA; 10% H |

*Figure 1. Summer ALEKS PPL scores for students enrolled in College Algebra (Fall term), 2016-18.*
*Gold: single assessment (SA); blue: best of two assessments; purple: best of more than two.*
*H: best score of 60 or greater (eligible for more advanced course); the vertical red line indicates this threshold.*

Roughly half of the UCSC students enrolling in Fall offerings of *Precalculus* enrollments are already calculus-eligible, or have already met the mathematics requirement for their intended major; this enrollment pattern has been strikingly persistent over several years. The ratios of calculus-eligible to appropriately placed *Precalculus* students are far worse for some demographic groups, particularly those vulnerable to stereotype threat. Among students who enrolled in *Precalculus* during or immediately after their 2018 Summer Orientation session, more than three quarters of the women, Hispanic/Latinx, and White/Causcasian students in majors requiring calculus who were already eligible to enroll in the calculus sequence, *Calculus with Applications*, required for their major (Table 3). Enrollment below placement was also widespread among students intending to major in Mathematics, Physics, or Engineering, but this behavior was far less common among these students.

|  | Life Sciences & Chemistry | Physics, Math & Eng. |
|---|---|---|
| Women | 78% | 29% |
| Hispanic/Latinx | 77% | 50% |
| White/Causcasian | 82% | 29% |
| Pell-eligible/First Gen | 64% | 37% |

*Table 3. Percentage of the students who enrolled in Precalculus during or shortly after 2018 Summer Orientation who were eligible to enroll in the calculus sequence required for their major.*

To look on the bright side, the number of students who are investing time, effort, and tuition in preparatory courses unlikely to increase their academic success has dramatically decreased. Unfortunately, the students most likely to benefit from preparatory math instruction—students who place into preparatory courses after full engagement with the placement process—are still disadvantaged by sharing those classes with large numbers of overqualified students. Distorted instructor expectations and competition for fixed resources (e.g. access to the instructor or TA's) can increase the challenges for the courses' intended audience. Comparison to other students in their demographic group can further erode the math self-efficacy of students who place into a preparatory course if an overwhelming majority of the other students sharing that course and demographic group are already proficient with the course content. Educating instructional staff about aggregate enrollment patterns and supporting implementation of pedagogies that support development of self-efficacy and minimize competition between students (e.g., use of criterion-based grading schemes, rather than percentile-based grading) can reduce the negative consequences of imperfect placement.

*Acknowledgements*
This research was supported by the Office of the President of the University of California, and by the Office of the Dean of Physical and Biological Sciences, University of California Santa Cruz.

# Appendix A.  2016 UCSC math placement outcomes

*A.1  Prior mathematical preparation*

Disaggregation by self-reported last high school math course suggests that reassessment captures typical retained/regained knowledge and skills from high school coursework, and that many students who reassess may have estimated what they "should" be able to do, based on their prior coursework, and reassessed because their initial score fell below that threshold. When considering mean assessment scores for students sharing self-reported final high school courses (Figure A.1), the average differences between scores for students who completed a single assessment and the best scores of students who completed multiple assessments differed by less than 5 points for all courses except Advanced Algebra, while the average improvement on reassessment was nearly 20 points; among students whose last high school math course was either AP Calculus or AP Statistics, the single/best of multiple assessments means differed by a single point.

| Last pre-UCSC math course | Single/multiple assessment ratio | Single assessment | Best of multiple assessments | Improvement on reassessment |
|---|---|---|---|---|
| Advanced Algebra | 2/1 | 53* | 69** | 25 |
| Precalculus | 3/2 | 68** | 71† | 19 |
| Statistics | 1/1 | 71† | 67** | 19 |
| AP Statistics | 2/1 | 74† | 75‡ | 19 |
| Calculus | 7/3 | 79‡ | 76‡ | 17 |
| AP/IB Calculus | 7/3 | 78‡ | 79‡ | 19 |

*Table A.1. Mean placement scores vs last high school math course for Fall 2016 survey respondents. Enrollment eligibility caps: **College Algebra (0-59)**, ****Precalculus (60-69)**, †Calculus with Applications (70-74), ‡Calculus for Math, Physics, and Engineering (75+).*

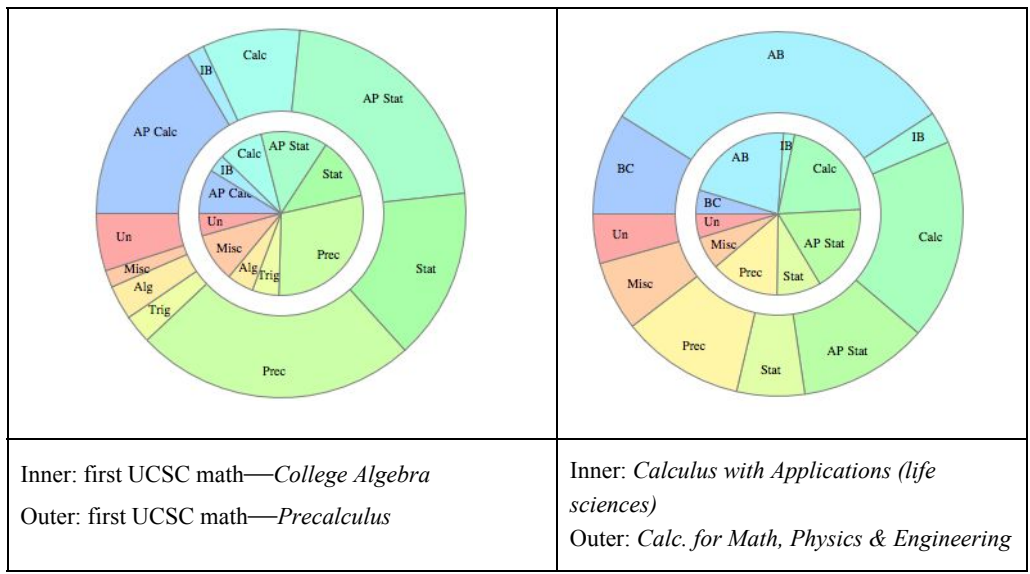

| Inner: first UCSC math—*College Algebra* | Inner: *Calculus with Applications (life sciences)* |
| Outer: first UCSC math—*Precalculus* | Outer: *Calc. for Math, Physics & Engineering* |

*Figure A.1. Last pre-UCSC math course, grouped by math enrollment in Fall 2016 (self-reported, survey). AB (BC): AP Calculus, AB (BC) exam; AP Calc: AP Calculus, exam unknown; IB: International Baccalaureate Higher Level Math; Prec: Precalculus; Alg: Advanced Algebra; Un: Unknown.*

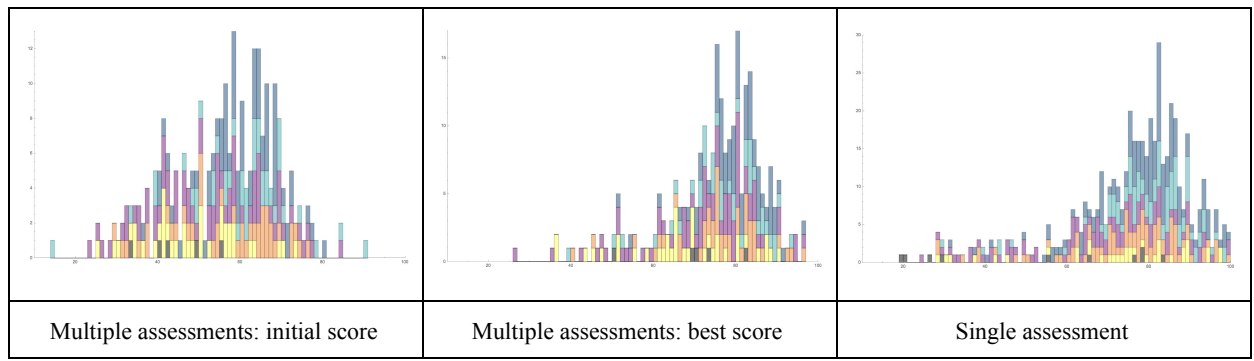

| Multiple assessments: initial score | Multiple assessments: best score | Single assessment |

*Figure A.2. Assessment scores for Fall 2016 survey respondents, by assessment activity and last high school course. Black: Advanced algebra, yellow: statistics, orange: AP stat., purple: precalculus, teal: calculus, blue: AP calculus.*

### A.2  Placement outcomes by ethnicity

Assessment outcomes for summer 2016, disaggregated by ethnicity, are summarized in Figure 2. Final placement tier distributions (single assessment or best of multiple assessments) for UCSC's white/Caucasian and Hispanic/Latinx students are relatively close, particularly in the middle three tiers; the percentage of students placing into the highest tier is higher for whites, while the percentage in the lowest tier is higher for Latinx students. Interestingly, the distributions of initial placements among students who reassessed at least once were nearly identical for white and Latinx students. Average assessment scores for Asian students are substantially higher than those for students in other ethnicity categories, regardless of the number of assessments.

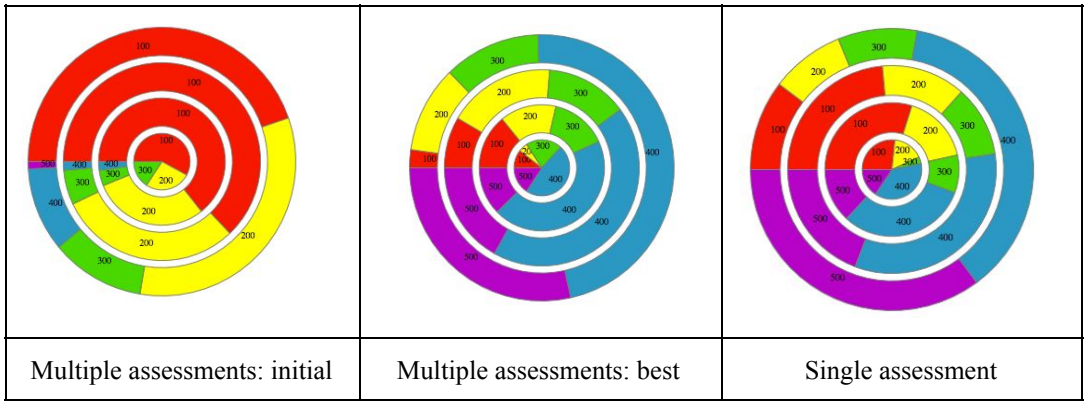

| Multiple assessments: initial | Multiple assessments: best | Single assessment |

*Figure A.3. Assessment outcomes by ethnicity, Summer 2016.*
*Outer ring: Asian; second ring: White/Caucasian; third ring: Hispanic/Latinx; inner pie: other.*
*Enrollment eligibility caps: 100 (red): College Algebra, 200 (yellow): Precalculus, 300 (green): Calculus with Applications, 400 (blue): Calculus for Math, Physics, and Eng., 500 (purple): Honors Calculus.*

Rates of reassessment, and of improvement on reassessment, were comparable across the main demographic groups. Reassessment provides a viable pathway to equity, narrowing placement gaps between ethnic groups: for example, the percentage of Latinx students who remained in the lowest placement tier after reassessment was only slightly larger than the percentage of Asian students who placed into that lowest tier and opted not to reassess.

## Appendix B. Starting point- and goal-oriented concise placement references

### ALEKS Score Placement

| ALEKS Score | MP Tier | Math Placement |
|---|---|---|
| Below 60 | 100 | MATH 2 |
| 60-69 | 200 | MATH 3 |
| 70-74 | 300 | MATH 11A |
| 75-84 | 400 | MATH 11A or MATH 19A |
| 85-100 | 500 | MATH 11A, 19A, or 20A |

*Figure B.1. Traditional placement conversion table*

# Getting Started in Math at UCSC

## Which first calculus course is right for you?

| Major | Calculus course | Course number | Minimum MP tier | Min score ALEKS | OR min AP score |
|---|---|---|---|---|---|
| 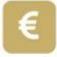 | Calculus for Economics | AMS 11A | 300 | 70 | Calculus AB 3 or BC 3 |
| 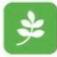 | Calculus with Applications | Math 11A | 300 | 70 | Calculus AB 3 or BC 3 |
| 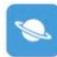 | Calculus for the Sciences, Mathematics & Engineering | Math 19A | 400 | 75 | Calculus AB 3 or BC 3 |
| 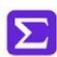 | Honors Calculus (not required for any major) | Math 20A | 500 | 85 | AB 4 or BC 3 |

## Psychology or Environmental Studies major?

The **General Psychology**, **Intensive Psychology**, and **Environmental Studies** majors math requirements are satisfied by:
- **MP tier of 300** or higher (ALEKS score of 70-100)    **OR**
- **AP Calculus** AB or BC score of 3 or higher    **OR**
- **AMS 3**, or Math 3 or higher math course.

**Cognitive Science** majors must complete one of the following calculus courses: AMS 11A, Math 11A, Math 19A, Math 20A.

## Not yet eligible for the course you need?

**Improve your Math Placement tier**: work in your ALEKS PPL Learning Module for at least 5 hours, and then **take another ALEKS placement assessment;**    **OR**

| To become eligible for | Pass preparatory course | Course number | Min MP for prep course | Min ALEKS score for prep course |
|---|---|---|---|---|
| AMS 11A or Math 11A | Precalculus for the Social Sciences | AMS 3 | 200 | 60 |
| AMS 11A or Math 11A or 19A | Precalculus | Math 3 | 200 | 60 |
| AMS 3 or Math 3 | College Algebra | Math 2 | 100 | 0 |

*Figure B.2. Goal-oriented choice architecture in one-page placement reference.*

22